\documentclass[a4]{amsart}

\usepackage{amssymb}
\usepackage[all]{xy}
\usepackage{amsmath, amsfonts, amsthm , amssymb} 
\usepackage{mathrsfs}
\usepackage{color}
\usepackage{hyperref}

\input xypic 
\usepackage{tikz} 
\usepackage{tikz-cd}
\usepackage{tikz-3dplot}

\usepackage{romannum}

\newtheorem{theorem}{Theorem}[section]
\newtheorem{proposition}[theorem]{Proposition}
\newtheorem{lemma}[theorem]{Lemma}
\newtheorem{corollary}[theorem]{Corollary}

\theoremstyle{definition}

\newtheorem{remark}[theorem]{Remark}

\newtheorem{question}[theorem]{Question}

\newcommand{\sk}{\ensuremath{\mbox{skel}\,}}

\newcounter{bean}

\newcommand{\namedright}[3]{\ensuremath{#1\stackrel{#2}
 {\longrightarrow}#3}}
\newcommand{\nameddright}[5]{\ensuremath{#1\stackrel{#2}
 {\longrightarrow}#3\stackrel{#4}{\longrightarrow}#5}}

\newcommand{\larrow}{\relbar\!\!\relbar\!\!\rightarrow}

\newcommand{\lnamedright}[3]{\ensuremath{#1\stackrel{#2}
 {\larrow}#3}}

\newcommand{\qqed}{\hfill\Box}

\newcommand{\Q}{\mathbb{Q}}

\newcommand{\Z}{\mathbb{Z}}

\newcommand{\dimp}{n}
\newcommand{\num}{\ell}

\begin{document}
\title[Homotopy rigidity]{Homotopy rigidity for quasitoric manifolds over a product of $d$-simplices} 

\author{Xin Fu}
\thanks{X.\,Fu is supported by the Beijing Natural Science Foundation (No.\,1244043).}
\address{Beijing Institute of Mathematical Sciences and Applications, Beijing 101408, China} 
\email{x.fu@bimsa.cn} 
\author{Tseleung So} 
\thanks{T.\,So is supported by NSERC Discovery Grant and NSERC RGPIN-2020-06428.}
\address{Department of Mathematics, University of Western Ontario, London ON, N6A 5B7, Canada} 
\email{tso28@uwo.ca} 
\author{Jongbaek Song} 
\address{Department of Mathematics Education, Pusan National University, Pusan 46241, Republic of Korea} 
\email{jongbaek.song@pusan.ac.kr} 
\author{Stephen Theriault}
\address{Mathematical Sciences, University of Southampton, Southampton SO17 1BJ, UK} 
\email{S.D.Theriault@soton.ac.uk} 


\begin{abstract} 
For a fixed integer $d\geq 1$, we show that two quasitoric manifolds over a product of $d$-simplices are homotopy equivalent after appropriate localization, provided that their integral cohomology rings are isomorphic.

\end{abstract}

\maketitle 

\pagenumbering{arabic}

\section{Introduction} 
\label{sec:intro}

A quasitoric manifold $M$ is a smooth, compact $2\dimp$-dimensional manifold endowed with a locally standard $T^{\dimp}$-action, such that the orbit space $M/T^{\dimp}$ is an $n$-dimensional simple polytope $P$. Here $T^{\dimp}=(S^1)^{\dimp}$ is the compact torus of rank ${\dimp}$. The Cohomological Rigidity Problem in toric topology~\cite{MS} poses the question of whether two quasitoric manifolds are homeomorphic or diffeomorphic if their integral cohomology rings are isomorphic. Special cases have provided positive evidence, such as four-dimensional quasitoric manifolds~\cite{OrRa, Fre}, Bott manifolds~\cite{CHJ}, certain generalized Bott manifolds~\cite{CMS-tr}, quasitoric manifolds with second Betti number equal to $2$~\cite{CPS}, and 6-dimensional quasitoric manifolds associated with Pogorelov polytopes~\cite{BEMPP} or a $3$-cube~\cite{Ha}. 

As cohomology rings are homotopy invariant, it is more natural to ask whether two quasitoric manifolds are homotopy equivalent if their cohomology rings are isomorphic. In theory, this should be a more accessible problem, while providing a test as to whether two quasitoric manifolds with isomorphic cohomology rings are homeomorphic or diffeomorphic.

\begin{question}
Let $M$ and $N$ be two quasitoric manifolds. If their integral cohomology rings are isomorphic, are they homotopy equivalent?
\end{question}

In fact, one can ask the question above for a broader class of spaces with torus actions, such as toric orbifolds \cite[Section~7]{DJ}. No counterexamples are known while several affirmative results have been established. For instance, one can deduce a positive answer for weighted projective spaces from the work in \cite{BFNR}, and the first three authors proved an affirmative answer for four dimensional toric orbifolds whose homology groups have no 2-torsion \cite{FSS}.



In this paper we focus on $2{\dimp}$-dimensional quasitoric manifolds $M$ with orbit space $P=\prod^{\num}_{i=1} \Delta^d$, a product of~$d$-dimensional simplices~$\Delta^d$, where $\dimp=\num d$. 
Notably, the class of quasitoric manifolds in this context includes certain generalised Bott manifolds~$B_{\num}$ for $d\geq 1$ that arise from a sequence 
\[
B_{\num}\xrightarrow{\pi_{\num}} B_{\num-1}\to \cdots \to B_1\xrightarrow{\pi_1} B_0= \text{ point},
\]
where each $\pi_i\colon B_i\to B_{i-1}$ for $i=1,\ldots, \num$ is a $\mathbb{C}P^d$-fibration.
However, not every quasitoric manifold over $\prod^{\num}_{i=1} \Delta^d$ is a generalized Bott manifold; see \cite[Section~5]{Ha}. 

A simple polytope $P$ can be associated with a moment-angle manifold~$\mathcal Z_P$. This manifold comes with a $T^m$-action, where $m$ is the number of facets (that is codimension-one faces)  of the polytope~$P$. A quasitoric manifold $M$ over $P$ can be regarded as a quotient
\[
M=\mathcal Z_P/T^{m-{\dimp}}
\]
by a freely acting subtorus $T^{m-{\dimp}}$ of $T^m$ on $\mathcal Z_P$, where~${\dimp}$ is the dimension of~$P$. 
This process results in a principal $T^{m-{\dimp}}$-fibration
\[
\nameddright{T^{m-\dimp}}{}{\mathcal{Z}_{P}}{}{M}.
\]


Throughout the paper, $H^*(M)$ denotes the integral cohomology ring of $M$ unless specified otherwise.
Now, we introduce our main result as follows.

\begin{theorem} 
   \label{main} 
   Let $M$ and $N$ be $2n$-dimensional quasi-toric manifolds with orbit space $\prod^{\num}_{i=1}\Delta^d$ for some $d\geq 1$ and let $\mathcal{P}$ be the set of primes $p\leq\dimp-d+1$. If there is a ring isomorphism $H^{\ast}(M)\cong H^{\ast}(N)$ then, after localizing away from~$\mathcal{P}$, there is a homotopy equivalence $M\simeq N$. 
\end{theorem} 


Two remarks should be made.
First, if $\Phi\colon M\to N$ is a homotopy equivalence obtained from Theorem~\ref{main}, then it induces an isomorphism $\Phi^*\colon H^*(N;\Q)\to H^*(M;\Q)$ and hence is a rational homotopy equivalence.
For specific values of $d$ and~$\ell$, Theorem~\ref{main} shows that this rational equivalence occurs after localizing away from a small number of primes.


Second, Theorem~\ref{main} should be compared to the original cohomological rigidity results for Bott manifolds~\cite{CHJ}, 2-stage generalized Bott manifolds~\cite{CMS-tr}, and quasitoric manifolds over a cube~\cite{Ha}. While we impose stronger conditions by localizing away from certain primes and consider only homotopy equivalences, the benefit of Theorem~\ref{main} is that it works for a larger class of quasitoric manifolds. For instance, if $d=1$ then $P=I^{\ell}$ is an $\ell$-dimensional hypercube. Two quasitoric manifolds $M$ and~$N$ over $I^{\ell}$ are homotopy equivalent if there is a ring isomorphism $H^*(M)\cong H^*(N)$ and localization occurs away from primes $p\leq\ell$. 
In particular, this works if $M$ and $N$ are Bott manifolds. 
Moreover, setting $d=1$ and $\ell=3$,  Theorem~\ref{main} gives a homotopy version of Hasui's result~\cite{Ha} after localizing away from $2$ and $3$.

The main result also works for $\ell$-stage generalized Bott manifolds $B_{\ell}$ that are constructed from iterated $\mathbb{C}P^d$-fibrations $\pi_i\colon B_i\to B_{i-1}$ for $i=1, \dots, \ell$ starting from $B_0=\{pt\}$. In this case, the associated simple polytopes are $P=\prod^{\ell}_{i=1}\Delta^d$. For instance, two $3$-stage generalized Bott manifolds over~$\prod_{i=1}^3\Delta^2$ having isomorphic cohomology rings are homotopy equivalent after localizing away from 2, 3 and~5. 
It is worth noting that the original cohomological rigidity problem in toric topology is true for 2-stage generalized Bott manifolds \cite[Theorem 1.3]{CMS-tr}, which are quasitoric manifolds over the product of two simplices. Theorem \ref{main} provides positive evidence that this extends to higher-stage generalized Bott manifolds.

\section{Preliminary information} 
\label{sec:prelim} 

\subsection{A review of quasitoric manifolds}

This section defines terms and makes some preliminary observations. We start with moment-angle manifolds from simple polytopes following Buchstaber and Panov~\cite[Section~6.2]{BP}. See also~\cite[Section 4]{DJ}.
Let $P$ be an $\dimp$-dimensional simple polytope. In other words  $P$ is a convex polytope having exactly~$\dimp$ facets intersecting at each vertex.  
Let 
\[
\mathcal F(P)=\{F_1,\ldots, F_m\}
\]
be the set of facets of $P$.

The moment-angle manifold $\mathcal Z_P$ is the quotient space
\[
\mathcal Z_P=P\times T^m/_\sim.
\]
Here $(x,t)\sim (x',t')$ if and only if $x=x'$ and $t^{-1}t'\in \prod_{i\in\mathcal I(x)} S^1_i$, where $\mathcal I(x)=\{i\mid x\in F_i\}$. 
There exists a~$T^m$-action on $\mathcal Z_P$ given by
\begin{equation}\label{torus action on zp}
T^m\times \mathcal Z_P\to \mathcal Z_P, \quad (g,[x,t])\mapsto [x,gt]
\end{equation}
for $g\in T^m$ and the  equivalence class $[x,t]\in \mathcal Z_{P}$ of $(x,t)\in P\times T^m$. 


Next, we define a quasitoric manifold following Davis and Januszkiewicz~\cite{DJ}. 
A $2\dimp$-dimensional manifold has a \emph{locally standard} $T^{\dimp}$-action if locally it is 
the standard action of $T^{\dimp}$ on $\mathbb{C}^{\dimp}$. A \emph{quasitoric manifold} 
over $P$ is a closed, smooth $2\dimp$-dimensional manifold $M$ that has a smooth locally standard 
$T^{\dimp}$-action for which the orbit space $M/T^{\dimp}$ is homeomorphic to $P$ as a manifold with corners.

A \emph{characteristic pair} $(P,\lambda)$ consists of an $\dimp$-dimensional simple polytope $P$ and a  function 
\[\lambda\colon \mathcal F(P)\to \Z^{\dimp}\]
 satisfying:
\begin{itemize}
\item $\lambda(F_i)$ is primitive for $i=1,\ldots,m$;
\item the set $\{\lambda(F_{i_1}),\ldots, \lambda(F_{i_k})\}$ extends to a basis of $ \Z^{\dimp}$ whenever $ F_{i_1}\cap \cdots \cap F_{i_k} \neq \emptyset$.
\end{itemize}
Such a function is called a \emph{characteristic function}. For a face $F=F_{i_1}\cap \cdots \cap F_{i_k}$ of codimension $k$, let $T_F$ denote the $k$-dimensional subtorus of $T^{\dimp}$ spanned by $\{\lambda(F_{i_1}),\ldots, \lambda(F_{i_k})\}$. If $F=P$, then~$T_F$ is the trivial subgroup.

For a characteristic pair $(P,\lambda)$, define the quotient space 
\[
M(P,\lambda)=P\times T^{\dimp}/_{\sim_\lambda}
\]
by the equivalence relation: $(x,t)\sim_{\lambda} (x',t')$ if and only if $x=x'$ and
$t^{-1}t'\in T_F$, where $F$ is the unique face such that $x=x'$ lies in its relative interior.
Notice that every quasitoric manifold~$M$ over $P$ can be constructed as such a quotient space. Here, $T^{\dimp}$ acts on the torus factor of~$P\times T^{\dimp}/_{\sim_\lambda}$, similarly to~\eqref{torus action on zp}.

Moment-angle manifolds and quasitoric manifolds are linked.
A characteristic function 
\[\lambda\colon \mathcal F(P)\to \Z^{\dimp}, \quad F_i\mapsto (\lambda_{1i},\ldots, \lambda_{\dimp i}), \] defines a linear map of lattices 
\[
\Lambda\colon \Z^m\to \Z^{\dimp}, \quad e_i\mapsto  (\lambda_{1i},\ldots, \lambda_{\dimp i}),
\]
where $\{e_1,\ldots, e_m\}$ is the standard basis of $\Z^m$.
Take the exponential of $\Lambda$ to get a homomorphism $\exp\Lambda\colon T^m\to T^{\dimp}$ of tori
sending $(t_1,\ldots, t_m)$ to $(t_1^{\lambda_{11}}t_2^{\lambda_{12}}\cdots t_m^{\lambda_{1m}},\ldots, t_1^{\lambda_{\dimp 1}}t_2^{\lambda_{\dimp 2}}\cdots t_m^{\lambda_{\dimp m}})$.
By~\cite[Proposition 7.3.13]{BP}, the kernel of $\exp\Lambda$ is isomorphic to $T^{m-n}$, which acts freely on $\mathcal Z_P$ and a $2n$-dimensional quasitoric manifold $M$ is $T^n$-equivariantly homeomorphic to the quotient $\mathcal{Z}_{P}/\ker\exp\Lambda$ equipped with the residual $T^n$-action. This implies that there is a principal $T^{m-\dimp}$-fibration 
\[\nameddright{T^{m-\dimp}}{}{\mathcal{Z}_{P}}{}{M}.\]  
As this is principal, it is classified by a map $M \to BT^{m-n}$
and there is a homotopy fibration 
\begin{equation} 
  \label{ZMfib} 
  \nameddright{\mathcal{Z}_{P}}{}{M}{}{BT^{m-\dimp}}.  
\end{equation} 

%

The cohomology of $M$ was calculated in~\cite{DJ}. They showed that there is a ring isomorphism 
\[
H^{\ast}(M)\cong\mathbb{Z}[x_{1},\ldots,x_{m}]/\mathcal I+\mathcal J
\] 
where each $x_{i}$ has degree $2$,  
the ideal $\mathcal I$ 
is generated by monomials 
$x_{i_1}\cdots x_{i_k}$ for which the intersection of $F_{i_1},\ldots,  F_{i_k}$ is empty, 
which is often called the \emph{Stanley--Reisner ideal} of $P$, 
 and $\mathcal J$ is an ideal of linear relations 
$\lambda_{j1}x_1+\lambda_{j2}x_2+\cdots+\lambda_{jm}x_m$ for  $j=1,\ldots n$. 
The cohomological properties of $M$ that will 
be relevant to us are: 
\begin{itemize} 
   \item $H^{2}(M)$ has rank $m-\dimp$; 
   \item $H^{\ast}(M)$ is multiplicatively generated by degree-two elements.
\end{itemize} 
Note that these two properties also imply that $M$ is simply-connected. We record a 
simple property of $M$ that follows immediately from the homotopy fibration~(\ref{ZMfib}) 
and the fact that $T^{m-\dimp}$ is an Eilenberg--MacLane space $K(\mathbb{Z}^{m-\dimp},1)$.

\begin{lemma} 
   \label{loopM} 
   Let $M$ be a quasitoric manifold associated to a characteristic pair $(P,\lambda)$. 
   Then the map $\mathcal{Z}_P \to M$
   in~(\ref{ZMfib}) induces an isomorphism $\pi_{t}(\mathcal{Z}_{P})\cong\pi_{t}(M)$ for $t\geq 3$.~$\qqed$ 
\end{lemma} 

We will need an identification of $\mathcal{Z}_{P}$ in the special case when $P=P_1\times P_2$ is a product of two simple polytopes $P_1$ and $P_2$.  The following statement can be found 
in~\cite[Proposition 6.4]{BP-old}; we give a brief proof. 

\begin{lemma} 
   \label{zkntype} 
   There is a homeomorphism 
   \[\mathcal Z_{P_1\times P_2}\cong\mathcal Z_{P_1}\times\mathcal Z_{P_2}.\]  
   In particular, if $P=\prod_{i=1}^{\num} \Delta^d$ is a product of $d$-simplices $\Delta^d$, then $\mathcal Z_{P}\cong\prod_{i=1}^{\num} S^{2d+1}$.
\end{lemma} 

\begin{proof} 
The first homeomorphism follows directly from the definition of a moment-angle manifold.
In the simplex case we have
$\mathcal Z_{\Delta^d}\cong S^{2d+1}$.
Therefore $\mathcal Z_{\prod_{i=1}^{\num} \Delta^d}\cong\prod_{i=1}^{\num} S^{2d+1}$.
\end{proof}

\subsection{Homotopy theory of quasitoric manifolds}

Consider the quasitoric manifold $M$ with orbit space $P=\prod_{i=1}^\ell \Delta^d$ as in Theorem~\ref{main}.
As~$\mathcal{Z}_{P}\cong\prod^{\num}_{i=1} S^{2d+1}$,  by~\eqref{ZMfib}, there is a homotopy fibration
\begin{equation}\label{eqn_M define fib}
\prod^{\num}_{i=1} S^{2d+1}\overset{}{\longrightarrow}M\overset{\delta}{\longrightarrow}BT^{\num}.
\end{equation}
Since $M$ is simply-connected, it has a $CW$-structure in which each cell corresponds to a homology class.
Fix such a $CW$-structure and for~$1\leq t\leq 2\dimp$, let $\sk_{t}(M)$ 
be the $t$-skeleton of~$M$.  
Since the cohomology of $M$ is concentrated in 
even degrees, there are homotopy equivalences \[
\sk_{2k+1}(M)\simeq\sk_{2k}(M)
\]
for any $k$. Therefore, we may focus only on the skeletons indexed by even integers. Consider the homotopy cofibrations 
\begin{equation}\label{eqn_M skeleton cofib}
\nameddright{\bigvee S^{2k-1}}{f_{k}}{\sk_{2k-2}(M)}{}{\sk_{2k}(M)}
\end{equation}
for $2\leq k\leq \dimp$. Note that $\sk_{0}(M)\simeq\ast$ and $\sk_{2\dimp}(M){\simeq} M$.

First we study the homotopy groups of $\sk_{2k}(M)$.

\begin{lemma}\label{lemma_hmtpy gp of M localized}
Let $M$ be a quasitoric manifold as in Theorem~\ref{main}, and let $\mathcal{P}$ be the set of primes $p\leq\dimp-d+1$.  Fix an integer $k$ such that $k\geq d+1$. If $2d+1<i\leq 2k$ then the homotopy group $\pi_i(\sk_{2k}(M))$ is a finite group of an order divisible only by primes in $\mathcal{P}$. 
Consequently, after localizing away from $\mathcal{P}$, there is an isomorphism  
\[
\pi_i(\sk_{2k}(M))\cong 0.
\]
\end{lemma}

\begin{proof}
The fibration sequence~\eqref{eqn_M define fib} induces a long exact sequence of homotopy groups. Since $\pi_i(BT^{\ell})$ is trivial for $i>2$, we have
\[
\pi_i(M)\cong\pi_i(\prod^{\num}_{j=1}S^{2d+1})\cong\bigoplus^{\num}_{j=1}\pi_i(S^{2d+1}).
\]
It is known that $\pi_i(S^{2d+1})$ is a finite group when $i\neq 2d+1$, and the first nontrivial $p$-torsion element in $\pi_*(S^{2d+1})$ occurs in $\pi_{2(p+d-1)}(S^{2d+1})$ (see~\cite{To}). By definition, $\mathcal{P}$ consists of all the primes that divide the order of $\pi_{t}(S^{2d+1})$ for any~$2d+1<t\leq 2\dimp$.
Hence the lemma holds true for $\sk_{2n}(M)\simeq M$.

Let $k$ be such that $d+1\leq k< n$. Let $F_{k}$ be the homotopy fibre of the skeletal inclusion 
\(\sk_{2k}(M)\to M\). 
Since the first cells in $M$ that are not in $\sk_{2k}(M)$ occur in dimension $2k+2$, the space $F_{k}$ is  
$2k$-connected. Therefore, 
\(\sk_{2k}(M)\to M\) 
induces an isomorphism $\pi_{i}(\sk_{2k}(M))\cong\pi_{i}(M)$ for all $i\leq 2k$. In particular, 
for $2d+1<i\leq 2k$, as $\pi_{i}(M)$ is a finite group of an order divisible only by primes in $\mathcal{P}$, 
the same is true for $\pi_{i}(\sk_{2k}(M))$.
\end{proof}

Next we record two cohomological properties for the map 
\(\delta\colon M\to BT^{\num}\).  

\begin{lemma}\label{lemma_gamma}
The ring homomorphism $\delta^*\colon H^{*}(BT^{\num})\to H^{*}(M)$ is a surjection 
that induces an isomorphism in degree $2$.
\end{lemma} 

\begin{proof} 
As $BT^{\num}\simeq\prod_{i=1}^{\num}\mathbb{C}P^{\infty}$, it follows that $H^{\ast}(BT^{\num})$ is a polynomial algebra 
generated by degree~$2$ elements. By~\eqref{eqn_M define fib}, the homotopy fibre of $\delta$ is at least $2$-connected. Therefore $\delta$ induces an isomorphism on $\pi_{2}$. Since $M$ and $BT^{\ell}$ are both simply-connected, the Hurewicz Theorem implies that $\delta$ induces an isomorphism on $H_{2}$. 
Simple-connectivity also implies, by the Universal Coefficient Theorem, that $H^{2}$ is the dual of $H_{2}$. 
Therefore $\delta^*\colon H^2(BT^{\num})\to H^2(M)$ is an isomorphism. Since $H^*(M)$ is generated as an algebra by degree $2$ elements, this implies that $\delta^*$ is surjective. 
\end{proof}

\begin{lemma}\label{lemma_skel of M = BT^k}
Let $\delta_{2d}\colon\sk_{2d}(M)\to\sk_{2d}(BT^{\num})$ 
be the restriction of 
\(\namedright{M}{\delta}{BT^{\num}}\) 
to $2d$-skeletons. Then it is a homotopy equivalence and induces an isomorphism in cohomology. 
\end{lemma} 

\begin{proof}  
In fibration sequence~\eqref{eqn_M define fib}, since the base is $1$-connected and the fibre is $2d$-connected, the Serre exact sequence 
(see~\cite[Theorem 6.4.4]{A}, for example) implies that $\delta$ induces an isomorphism in 
homology in degrees~$\leq 2d$. Thus $\delta_{2d}$ induces an isomorphism in homology 
in all degrees, implying that it is a homotopy equivalence by Whitehead's theorem.  
\end{proof}

\section{The strategy for the proof of Theorem~\ref{main}} 
\label{sec:strategy} 
For $d\geq 1$, let $\Delta^{d}$ be a $d$-simplex and let $P=\prod^{\num}_{i=1}\Delta^{d}$ be a product of $d$-simplices. 
For two quasitoric manifolds $M$ and $N$ over $P$, assume that there is a ring isomorphism
\begin{equation}\label{eqn_gamma ring isom}
\gamma\colon H^{\ast}(N;\Z)\to H^{\ast}(M;\Z).
\end{equation}

The goal is to determine whether there is a homotopy equivalence $M\simeq N$. 
We will show that this is true after localizing away from $\mathcal{P}$. A similar 
local approach was taken in~\cite{Th} in the special case of Bott manifolds when $d=1$; our methods are different and generalize the result to all quasitoric manifolds over products of $d$-simplices for $d\geq 1$.

The argument proceeds by an induction on skeletons.
Consider the homotopy cofibrations 
\[\nameddright{\bigvee S^{2k-1}}{f_{k}}{\sk_{2k-2}(M)}{}{\sk_{2k}(M)}\] 
\[\nameddright{\bigvee S^{2k-1}}{f'_{k}}{\sk_{2k-2}(N)}{}{\sk_{2k}(N)}\] 
for $2\leq k\leq \dimp$.
We will construct a homotopy equivalence
\[
\Phi_{2k}\colon \sk_{2k}(M)\longrightarrow\sk_{2k}(N)
\]
by different methods in two cases: (\romannum{1}) $k=d$ and (\romannum{2}) $d+1\leq k\leq n$.
\medskip 

\noindent 
\textit{Case 1: $k=d$}. 
In {Proposition}~\ref{Phi_2d}, it will be shown that 
the ring isomorphism $\gamma\colon H^{\ast}(M)\to H^{\ast}(N)$ 
gives rise to a homotopy equivalence $\Phi_{2d}\colon\sk_{2d}(M)\longrightarrow\sk_{2d}(N)$. 
No localization is needed. 
\medskip 

\noindent 
\textit{Case 2: $d+1\leq k\leq n$}. 
In Proposition~\ref{k>2} an induction on skeletons  and novel methods relating the cohomology and homotopy groups of $CW$-complexes with only even dimensional cells 
will be used to show that, after localization away from~$\mathcal{P}$, there is a homotopy equivalence 
$\Phi_{2k}\colon\sk_{2k}(M)\to\sk_{2k}(N)$. 

\medskip 
 
Granting Propositions~\ref{Phi_2d} and~\ref{k>2}, we can prove the main result in the paper. 

\begin{proof}[Proof of Theorem~\ref{main}] 
Since $\sk_{2n}(M)\simeq M$ and $\sk_{2n}(N)\simeq N$, the $k=\dimp$ instance of \textit{Case 2} implies that after localizing away from $\mathcal{P}$ there is a homotopy 
equivalence $M\simeq N$. 
\end{proof}

\section{Construction of the homotopy equivalence \texorpdfstring{$\Phi_{2d}$}{}}\label{sec:k=d}

Given a ring isomorphism $\gamma\colon H^*(N)\to H^*(M)$, its restriction to degree $2$ is 
a group isomorphism  
\[
\gamma\colon H^2(N)\cong\Z^{\num}\to H^2(M)\cong\Z^{\num}.
\]
Apply the classifying space functor $B(-)$ twice to obtain a homotopy equivalence 
\begin{equation}\label{defn Gamma}
\Gamma\colon BT^{\num}\to BT^{\num}.
\end{equation}
Recall from \eqref{eqn_M define fib} that there are homotopy fibrations 
\[
\nameddright{\prod^{\num}_{i=1} S^{2d+1}}{}{M}{\delta}{BT^{\num}}
\quad
\text{and}\quad 
\prod^{\num}_{i=1} S^{2d+1}\overset{}{\longrightarrow}N\overset{\delta'}{\longrightarrow}BT^{\num}.
 \]

\begin{lemma}\label{lemma_delta^* commute}
There is a commutative diagram
\begin{equation}\label{diagram_gamma and delta^* commute}
\begin{split}
\xymatrix{
H^*(BT^{\num})\ar[r]^-{(\delta')^*}\ar[d]^-{\Gamma^*}	&H^*(N)\ar[d]^-{\gamma}\\
H^*(BT^{\num})\ar[r]^-{\delta^*}							&H^*(M).
}\end{split}
\end{equation}
\end{lemma}

\begin{proof} 
Observe that Diagram~\eqref{diagram_gamma and delta^* commute} commutes in degree $2$ by definition of $\Gamma$. The fact that both $H^{\ast}(BT^{\num})$ and $H^{\ast}(M)$ are generated 
as algebras by degree $2$ elements then implies the diagram commutes in all degrees since all maps are algebra maps. 
\end{proof}

Restricting 
\(\namedright{BT^{\num}}{\Gamma}{BT^{\num}}\) 
to $2d$-skeletons gives a map 
\[\Gamma_{2d}\colon \sk_{2d}(BT^{\num})\to\sk_{2d}(BT^{\num}).\] 
Using the homotopy equivalences 
\(\lnamedright{\sk_{2d}(M)}{\delta_{2d}}{\sk_{2d}(BT^{\num})}\)  
and 
\(\lnamedright{\sk_{2d}(N)}{\delta'_{2d}}{\sk_{2d}(BT^{\num})}\)  
from Lemma~\ref{lemma_skel of M = BT^k}, define the map $\Phi_{2d}$ by the composite 
\begin{equation}\label{eq_Phi_2d}
\Phi_{2d}\colon\sk_{2d}(M)\xrightarrow{\delta_{2d}}\sk_{2d}(BT^{\num}) 
   \xrightarrow{\Gamma_{2d}}\sk_{2d}(BT^{\num})\xrightarrow{(\delta'_{2d})^{-1}}\sk_{2d}(N).
\end{equation}
By definition of $\Phi_{2d}$, there is a commutative diagram 
\begin{equation}\label{diagram_gamma and delta^* commute 2d}
\begin{tikzcd}
 H^*(\sk_{2d}(BT^\ell)) \arrow[r,"(\delta'_{2d})^*"] \arrow[d,"\Gamma_{2d}^*"] 
 & H^*(\sk_{2d}(N)) \arrow[d,"\Phi_{2d}^*"]  \\
H^*(\sk_{2d}(BT^\ell)) \arrow[r,"\delta_{2d}^*"]  
 & H^*(\sk_{2d}(M)).
\end{tikzcd}
\end{equation} 

\begin{proposition}\label{Phi_2d} 
   The map 
   \(\namedright{\sk_{2d}(M)}{\Phi_{2d}}{\sk_{2d}(N)}\) 
   is a homotopy equivalence. 
\end{proposition} 

\begin{proof} 
Recall from~\eqref{defn Gamma} that $\Gamma$ is a homotopy equivalence, so $\Gamma^*$ is an isomorphism and 
so is~$\Gamma_{2d}^{\ast}$. As $H^{\ast}(BT^{\num})$ is concentrated in even degrees, 
the Universal Coefficient Theorem implies that $H_{\ast}(BT^{\ell})$ is dual 
to $H^{\ast}(BT^{\ell})$, and therefore  $(\Gamma_{2d})_{\ast}$ is also an isomorphism. 
Hence~$\Gamma_{2d}$ is a homotopy equivalence by Whitehead's Theorem. 
By Lemma~\ref{lemma_skel of M = BT^k} both $\delta_{2d}$ and $\delta'_{2d}$ are homotopy equivalences. Therefore the composite $\Phi_{2d}$ is a homotopy equivalence. 
%
\end{proof}

\section{Construction of the homotopy equivalence \texorpdfstring{$\Phi_{2k}$}{} for \texorpdfstring{$k\geq d+1$}{}}
We begin with a general lemma that is of interest in its own right. 

\begin{lemma}\label{attach_detect} 
Let $X$ be a connected CW complex having only even dimensional cells. For each $k\geq 1$ there is a group homomorphism
\[
g_X\colon H_{2k+2}(X)\to\pi_{2k+1}(\sk_{2k}(X))
\]
sending the homology class of a $(2k+2)$-cell $[e_{\alpha}]$ to the homotopy class of its attaching map $f_{\alpha}\colon\partial e_{\alpha}\to\sk_{2k}(X)$. This satisfies the following properties: 
\begin{itemize}  
   \item[(a)] if the standard maps of pairs 
            $H_{2k+2}(X)\to H_{2k+2}(X,\sk_{2k}(X))$
            and 
            $\pi_{2k+2}(X,\sk_{2k}(X))\to \pi_{2k+1}(\sk_{2k}(X))$  
            are both isomorphisms then so is $g_{X}$; 
   \item[(b)] if $Y$ is another CW-complex having only even dimensional 
            cells and $f\colon X\to  Y$ is a map then there is a commutative diagram 
            \begin{equation}\label{diagram_morphism g naturality}
              \begin{split}  \xymatrix{
                 H_{2k+2}(X)\ar[r]^-{g_X}\ar[d]^-{f_*}	&\pi_{2k+1}(\sk_{2k}(X))\ar[d]^-{(f_{2k})_*}\\
                 H_{2k+2}(Y)\ar[r]^-{g_Y}				&\pi_{2k+1}(\sk_{2k}(Y)).
                 }\end{split}
               \end{equation}
\end{itemize} 
\end{lemma}

\begin{proof}
First we define $g_X$.
The pair $(X, \sk_{2k}(X))$ induces a long exact sequence of relative homotopy groups
\[
\cdots\longrightarrow\pi_j(\sk_{2k}(X))\longrightarrow\pi_j(X)\overset{\jmath_{\pi}}{\longrightarrow}\pi_{j}(X,\sk_{2k}(X))\overset{\partial_{\pi}}{\longrightarrow}\pi_{j-1}(\sk_{2k}(X))\longrightarrow\cdots
\]
and a long exact sequence of relative homology groups
\[
\cdots\longrightarrow H_j(\sk_{2k}(X))\longrightarrow H_j(X)\overset{\jmath_H}{\longrightarrow}H_{j}(X,\sk_{2k}(X))\overset{\partial_{H}}{\longrightarrow}H_{j-1}(\sk_{2k}(X))\longrightarrow\cdots
\]
Since $X$ has cells only in even dimensions, the pair $(X,\sk_{2k}(X))$ is $(2k+1)$-connected and $\sk_{2k}(X)$ 
is simply-connected. Therefore, by~\cite[Proposition~4.28]{H} for example, there is an isomorphism  
\[\pi_{2k+2}(X,\sk_{2k}(X))\cong\pi_{2k+2}(X/\sk_{2k}(X)).\]  
On the other hand, the usual isomorphism $H_{m}(X,A)\cong \widetilde{H}_{m}(X/A)$ for 
pairs of spaces $(X,A)$ and any $m\geq 0$ implies in our case that there is an isomorphism 
\[H_{2k+2}(X,\sk_{2k}(X))\cong H_{2k+2}(X/\sk_{2k}(X)).\] 
Observe that the inclusion $\bigvee S^{2k+2}\hookrightarrow X/\sk_{2k}(X)$ of the bottom cells is a $(2k+3)$-equivalence, implying that there are isomorphisms  
\[
\pi_{2k+2}(X/\sk_{2k}(X))\cong\pi_{2k+2}(\bigvee S^{2k+2})
\quad\text{and}\quad
H_{2k+2}(X/\sk_{2k}(X))\cong H_{2k+2}(\bigvee S^{2k+2}).
\]
The Hurewicz homomorphism $\pi_{2k+2}(X/\sk_{2k}(X))\to H_{2k+2}(X/\sk_{2k}(X))$
is therefore an isomorphism. Combining these isomorphisms then gives an isomorphism 
\[hur\colon\pi_{2k+2}(X,\sk_{2k}(X))\to H_{2k+2}(X,\sk_{2k}(X)).\] 
Now define~$g_{X} \colon H_{2k+2}(X)\to\pi_{2k+1}(\sk_{2k}(X))$ by the composite 
\[
g_X\colon H_{2k+2}(X)\overset{\jmath_H}{\longrightarrow}H_{2k+2}(X,\sk_{2k}(X))\xrightarrow{hur^{-1}}
\pi_{2k+2}(X,\sk_{2k}(X))\overset{\partial_{\pi}}{\longrightarrow}\pi_{2k+1}(\sk_{2k}(X)).
\] 
Since $g_{X}$ is the composite of three group homomorphisms it too is a homomorphism, and by construction~$g_{X}$ 
sends the homology class of a $(2k+2)$-cell $e_{\alpha}$ to its attaching map 
$f_{\alpha}\in\pi_{2k+1}(\sk_{2k}(X))$. 

For part~(a), the composite defining $g_{X}$ implies that if both $\jmath_{H}$ and $\partial_{\pi}$ 
are isomorphisms then, as $hur$ is also an isomorphism, so is $g_{X}$. 

For part~(b), consider the diagram
\[
\xymatrix{ 
H_{2k+2}(X)\ar[r]^-{\jmath_H}\ar[d]^-{f_*}	&H_{2k+2}(X,\sk_{2k}(X))\ar[r]^-{hur^{-1}}\ar[d]^-{f_*}  
  & \pi_{2k+2}(X,\sk_{2k}(X))\ar[r]^{\partial_{\pi}}\ar[d]^-{f_*}  &\pi_{2k+1}(\sk_{2k+2}(X))\ar[d]^-{(f_{2k+2})_*}\\ 
H_{2k+2}(Y)\ar[r]^-{\jmath_H} &H_{2k+2}(Y,\sk_{2k}(Y))\ar[r]^{hur^{-1}} 
 & \pi_{2k+2}(Y,\sk_{2k}(Y))\ar[r]^-{\partial_{\pi}} &\pi_{2k+1}(\sk_{2k+2}(Y)). 
}
\] 
The left, middle and right squares commute by the naturality of $\partial_{\pi},\jmath_H$ and $hur$ 
respectively. The composites along the top and bottom rows are the definitions of $g_{X}$ and $g_{Y}$ 
respectively. Thus we obtain the commutativity of~\eqref{diagram_morphism g naturality}. 
\end{proof}

We apply Lemma~\ref{attach_detect} to the case $M\xrightarrow{\delta} BT^{\num}$ in~\eqref{eqn_M define fib}. Recall that $\mathcal{P}$ is the set of primes $p\leq\dimp-d+1$.

\begin{corollary}\label{cor_g_M and g_BT}
For $k\geq 1$, there is a commutative diagram
\begin{equation}\label{diagram_morphism g and BGamma}
\begin{split}
\xymatrix{
H_{2k+2}(M)\ar[r]^-{g_M}\ar[d]^-{\delta_*}	&\pi_{2k+1}(\sk_{2k}(M))\ar[d]^-{(\delta_{2{k}})_*}\\
H_{2k+2}(BT^{\num})\ar[r]^-{g_{BT^{\num}}}			&\pi_{2k+1}(\sk_{2k}(BT^{\num})),
}
\end{split}
\end{equation}
where $g_M$ is an injection and $g_{BT^{\num}}$ is an isomorphism.
Furthermore, if $k>d$ then
\[
g_{M}\colon H_{2k+2}(M)\to\pi_{2k+1}(\sk_{2k}(M))
\]
is an isomorphism after localizing away from $\mathcal{P}$.
\end{corollary}
\begin{proof}
The commutativity of~(\ref{diagram_morphism g and BGamma}) is immediate from Lemma~\ref{attach_detect} by taking $f$ to 
be~$\delta\colon M\to BT^{\num}$. 

To show that $g_{BT^\ell}$ is an isomorphism, we consider the pair $(BT^{\num},\sk_{2k}(BT^{\num}))$, which induces an isomorphism 
$H_{2k+2}(BT^{\num}) \to H_{2k+2}(BT^{\num},\sk_{2k}(BT^{\num}))$ for skeletal reasons. It also induces a long exact sequence of relative homotopy groups
\[
\begin{tikzcd}[column sep=small]
\cdots \arrow{r} &
\pi_j(BT^{\num}) 
\arrow{r} & 
\pi_{j}(BT^{\num},\sk_{2k}(BT^{\num})) \arrow{r} &
\pi_{j-1}(\sk_{2k}(BT^{\num})) \arrow{r} &
\pi_{j-1}(BT^{\num})\arrow{r} & 
\cdots. 
\end{tikzcd}
\]
For $k\geq 1$, both 
$\pi_{2k+2}(BT^{\num})$ and $\pi_{2k+1}(BT^{\num})$ are trivial, which implies that 
 $\pi_{2k+2}(BT^{\num},\sk_{2k}(BT^{\num}))$ and $\pi_{2k+1}(\sk_{2k}(BT^{\num}))$ are isomorphic. 
Thus, the map $g_{BT^{\num}}$ is an isomorphism by Lemma~\ref{attach_detect}(a). 

Next, the map $\delta_*\colon H_{2{k}+2}(M)\to H_{2{k}+2}(BT^{\num})$ 
in~(\ref{diagram_morphism g and BGamma}) is an injection by Lemma~\ref{lemma_gamma}. 
Therefore, the commutativity of~(\ref{diagram_morphism g and BGamma}) implies that $g_M$ is an injection.

Finally, suppose that $k>d$. By Lemma~\ref{lemma_hmtpy gp of M localized}, after localizing away from $\mathcal{P}$, homotopy groups $\pi_{s}(M)$ are trivial for $2d+1<s\leq 2n$. As~$d<k<n$ we have $2d+1<2k+1\leq 2n$, so both $\pi_{2k+1}(M)$ and $\pi_{2k+2}(M)$ are trivial. Therefore the boundary map
$\pi_{2k+2}(M,\sk_{2k}(M))\to\pi_{2k+1}(\sk_{2k}(M))$ in the long exact sequence of relative homotopy groups is an isomorphism. Since $\sk_{2k}(M)$ has trivial homology in degrees larger than $2k$, the map 
$H_{2k+2}(M)\to H_{2k+2}(M,\sk_{2k}(M))$ in the long exact sequence of relative homology groups is also an isomorphism. Therefore, by Lemma~\ref{attach_detect}~(a), $g_M$ is an isomorphism.
\end{proof}

\begin{lemma}\label{pre_lemma_isom btw attaching maps}
If there is a ring isomorphism $\gamma\colon H^*(N)\to H^*(M)$
then there is a commutative diagram 
\begin{equation}\label{eq_diag_MN}
\begin{split}
\xymatrix{ 
    H_{2d+2}(M)\ar[r]^-{g_{M}}\ar[d]^{\gamma^{\vee}} 
         & \pi_{2d+1}(\sk_{2d}(M))\ar[d]^{(\Phi_{2d})_{\ast}} \\
    H_{2d+2}(N)\ar[r]^-{g_{N}} & \pi_{2d+1}(\sk_{2d}(N))
}\end{split}
\end{equation}
where $\gamma^{\vee}$ is  the dual of $\gamma$.
\end{lemma} 
\begin{proof} 
As in~(\ref{defn Gamma}), $\gamma$ induces a map 
\(\Gamma\colon BT^{\ell}\to BT^{\ell}\). 
Consider the diagram 
\begin{equation}\label{eq_cub_diag}
\begin{tikzcd}
H_{2d+2}(BT^\ell) \arrow{rrr}{g_{BT^\ell}} \arrow{ddd}[swap]{\Gamma_\ast} \arrow[phantom, xshift=-10, yshift=3]{dddr}{\footnotesize\text{(A)}}
\arrow[phantom, xshift=-25, yshift=5]{drrr}{\footnotesize{\text{(D)}}} 
&&& \pi_{2d+1}(\sk_{2d}(BT^\ell))\arrow{ddd}{(\Gamma_{2d})_\ast} 
\arrow[phantom, xshift=10, yshift=3]{dddl}{\footnotesize\text{(C)}} \\
~& H_{2d+2}(M) \arrow{r}{g_M} \arrow{d}{\gamma^\vee}\arrow{ul}{\delta_\ast} & \pi_{2d+1}(\sk_{2d}(M))\arrow{d}{(\Phi_{2d})_\ast} \arrow{ur}[swap]{(\delta_{2d})_\ast}& ~\\
& H_{2d+2}(N) \arrow{r}{g_N} \arrow{dl}[swap]{\delta'_\ast} & \pi_{2d+1}(\sk_{2d}(N))\arrow{dr}{(\delta'_{2d})_\ast}&~ \\
H_{2d+2}(BT^\ell) \arrow{rrr}[swap]{g_{BT^\ell}} 
\arrow[phantom, xshift=-25, yshift=-5]{urrr}{\footnotesize{\text{(B)}}} 
& ~ &~& \pi_{2d+1}(\sk_{2d}(BT^\ell)).
\end{tikzcd}
\end{equation}
Observe that Diagram \eqref{eq_diag_MN} is located in the inner square. The outer square is obtained by applying Lemma~\ref{attach_detect} to the map $\Gamma\colon BT^\ell \to BT^\ell$, and hence it is commutative. Diagram (A) is the dual of \eqref{diagram_gamma and delta^* commute}, and hence it is also homotopy commutative. Diagrams  (B) and  (D) are commutative by Corollary~\ref{cor_g_M and g_BT}.  Diagram (C) commutes by the definition of $\Phi_{2d}$, see \eqref{eq_Phi_2d}. Therefore, from the commutativity of~(\ref{eq_cub_diag}), we obtain  
\begin{align*}
(\delta'_{2d})_\ast \circ g_N \circ \gamma^\vee &=g_{BT^\ell} \circ \Gamma_\ast \circ \delta_\ast  =(\Gamma_{2d})_\ast \circ g_{BT^\ell} \circ \delta_\ast = (\delta'_{2d})_\ast  \circ (\Phi_{2d})_\ast \circ g_M. 
\end{align*}
Since $\delta'_{2d}$ is a homotopy equivalence by Lemma~\ref{lemma_skel of M = BT^k}, 
 we then obtain $g_N \circ \gamma^\vee=(\Phi_{2d})_\ast \circ g_M$, which is the equality asserted by the lemma.
\end{proof} 

To construct homotopy equivalences $\Phi_{2k}\colon \sk_{2k}(M) \to \sk_{2k}(N)$ for $k\geq d+1$, we prepare more explicit notations for the map $g_M\colon H_{2k+2}(M) \to \pi_{2k+1}(\sk_{2k}(M))$ of Lemma \ref{attach_detect}.  For each $k\geq  1$, enumerate $(2k+2)$-cells of $M$ as $e_1, \dots, e_{s_k}$, whose attaching maps are $f_1, \dots, f_{s_k}$, respectively. To each $e_{i}$ there is a homology class $[e_i] \in H_{2k+2}(M)$, giving a group isomorphism  
$H_{2k+2}(M)\cong\mathbb{Z}\langle [e_{1}], \ldots ,[e_{s_k}]\rangle$. Define a linear map 
\begin{equation}\label{lineariso} 
    \Z\langle{f_1,\ldots,f_{s_k}}\rangle\longrightarrow\Z\langle  [e_{1}], \ldots ,[e_{s_k}] \rangle \cong H_{2k+2}(M) 
\end{equation} 
by sending $f_{i}$ to $[e_{i}]$. Since $g_{M}$ sends $[e_{i}]$ to the homotopy class of its attaching map, 
the composite 
\begin{equation}\label{eq_map_g}
g\colon
\Z\langle{f_1,\ldots,f_{s_k}}\rangle\longrightarrow\Z\langle [e_{1}],\ldots,[e_{s_k}] \rangle\cong
H_{2k+2}(M)
\xrightarrow{g_M}
\pi_{2k+1}(\sk_{2k}(M))
\end{equation}
sends $f_{i}$ to its homotopy class. 

The same logic applies to $N$, so that we have a group isomorphism $H_{2k+2}(N) \cong \mathbb{Z}\left< [e'_1], \dots, [e'_{s_k}] \right>$ where $[e'_1], \dots, [e'_{s_k}]$ are homology classes of cycles representing $(2k+2)$-cells $e'_1, \dots, e'_{s_k}$ of $N$. The composite
\begin{equation}\label{eq_map_g_prime}
g'\colon\Z\langle{f'_1,\ldots,f'_{s_k}}\rangle\longrightarrow\Z\langle [e'_{1}],\ldots,[e'_{s_k}] \rangle\cong
H_{2k+2}(N)
\xrightarrow{g_N}
\pi_{2k+1}(\sk_{2k}(N))
\end{equation}
sends $f'_i$ to its homotopy class.  



We are now ready to construct homotopy equivalences.

\begin{lemma}\label{lemma_inductive argument}
Given $k\geq 1$, consider the maps $g$ and $g'$ defined above and suppose there is a homotopy equivalence $\Phi_{2k}\colon\sk_{2k}(M)\to\sk_{2k}(N)$.
If there is a linear isomorphism $A\colon\Z\langle{f_1,\ldots,f_{s_k}}\rangle\to\Z\langle{f'_1,\ldots,f'_{s_k}}\rangle$ making the diagram
\begin{equation}\label{diagrm_commutative diagram for inductive argument}
\begin{split}\xymatrix{
\Z\langle{f_1,\ldots,f_{s_k}}\rangle\ar[r]^-{g}\ar[d]^-{A}	&\pi_{2k+1}(\sk_{2k}(M))\ar[d]^-{(\Phi_{2k})_*}\\
\Z\langle{f'_1,\ldots,f'_{s_k}}\rangle\ar[r]^-{g'}			&\pi_{2k+1}(\sk_{2k}(N))
}\end{split}
\end{equation}
commute, then there is a homotopy equivalence $\Phi_{2k+2}\colon\sk_{2k+2}(M)\to\sk_{2k+2}(N)$.
\end{lemma}

\begin{proof} 
The linear isomorphism 
\(A\colon\Z\langle{f_1,\ldots,f_{s_k}}\rangle\to\Z\langle{f'_1,\ldots,f'_{s_k}}\rangle\)  
can be geometrically realized by a homotopy equivalence  
$\varphi_A\colon\bigvee_{i=1}^{s_k} S^{2k+1}\to\bigvee_{i=1}^{s_k}  S^{2k+1}$  
as follows. If $A(f_{i})=\sum_{j=1}^{s_k} a_{ij}f'_{j}$ for coefficients $a_{ij}\in\mathbb{Z}$, 
then define $\varphi_{A}$ on the $i^{th}$-wedge summand to be the sum $\Sigma_{j=1}^{s_k} a_{ij}\iota_{j}$ 
where $a_{ij}$ is a map of degree $a_{ij}$ and 
$\iota_{j}\colon S^{2k+1}\to \bigvee_{l=1}^{s_k} S^{2k+1}$ 
is the inclusion of the $j^{th}$-wedge summand. We claim  
that~(\ref{diagrm_commutative diagram for inductive argument}) is then geometrically realized by
a homotopy commutative diagram
\begin{equation}\label{geomreal} 
\begin{split}
\xymatrix{
\bigvee_{i=1}^{s_k} S^{2k+1}\ar[r]^-{\bigvee_{i=1}^{s_k} f_i}\ar[d]^-{\varphi_A}	&\sk_{2k}(M)\ar[d]^-{\Phi_{2k}}\\
\bigvee_{i=1}^{s_k} S^{2k+1}\ar[r]^-{\bigvee_{i=1}^{s_k} f_i'}					&\sk_{2k}(N).
}\end{split}
\end{equation} 
To see this, restrict to the $i^{th}$ sphere of $\bigvee^{s_k}_{i=1}S^{2k+1}$. By~\eqref{diagrm_commutative diagram for inductive argument}, we have 
$\left( (\Phi_{2k})_*\circ g\right) (f_{i})=\left( g'\circ A\right)(f_{i})$. Since $A(f_{i})=\sum_{j=1}^{s_k} a_{ij}f'_{j}$ 
and~$g'$ is linear, we obtain 
$\left((\Phi_{2k})_*\circ g\right) (f_{i})=\left( \sum_{j=1}^{s_k} a_{ij}g'\right) (f_{j})$. Since $g$ sends $f_{i}$ to its homotopy 
class and $g'$ sends each $f'_{j}$ to its homotopy class, we obtain 
$\Phi_{2k}\circ f_{i}\simeq\sum_{j=1}^{s_k} a_{ij} f'_{j}$. The right side may be rewritten as 
$\left(\bigvee^{s_k}_{j=1}f'_j\right)\circ(\sum^{s_k}_{j=1}a_{ij}\iota_j)$, which by definition of $\varphi_{A}$ 
then equals $\left(\bigvee^{s_k}_{j=1}f'_j\right)\circ\varphi_{A}\circ\iota_{i}$. Thus~\eqref{geomreal} 
homotopy commutes when restricted to the $i^{th}$-wedge summand. As~$i$ was arbitrary, 
\eqref{geomreal} homotopy commutes.  

The homotopy commutativity of~\eqref{geomreal} implies that there is a homotopy cofibration diagram 
\[
\xymatrix{
\bigvee_{i=1}^{s_k} S^{2k+1}\ar[r]^-{\bigvee_{i=1}^{s_k} f_i}\ar[d]^-{\varphi_A}	&\sk_{2k}(M)\ar[d]^-{\Phi_{2k}}\ar[r]	&\sk_{2k+2}(M)\ar[d]^-{\Phi_{2k+2}}\\
\bigvee_{i=1}^{s_k} S^{2k+1}\ar[r]^-{\bigvee_{i=1}^{s_k} f_i'}					&\sk_{2k}(N)\ar[r]						&\sk_{2k+2}(N)
}
\]
where $\Phi_{2k+2}$ is an induced map of cofibres.
Since $\varphi_A$ and $\Phi_{2k}$ induce isomorphisms in homology, so does $\Phi_{2k+2}$ 
by the Five Lemma. Since all spaces are simply-connected, $\Phi_{2k+2}$ is therefore a homotopy equivalence by Whitehead's Theorem.
\end{proof} 

\begin{remark}\label{Aremark}  
There is a local version of Lemma~\ref{lemma_inductive argument}: in the statement 
and proof simply localize spaces away from $\mathcal{P}$ and change all instances of 
$\mathbb{Z}$ to the integers localized away from $\mathcal{P}$.
\end{remark}

\begin{proposition}\label{k>2}
Let $M$ and $N$ be $2n$-dimensional quasitoric manifolds with orbit space 
$P=\prod_{i=1}^{\num}\Delta^d$. If there is a ring isomorphism $H^{\ast}(M)\cong H^{\ast}(N)$, 
then after localizing away from $\mathcal{P}$, there is a homotopy equivalence 
$\Phi_{2k}\colon\sk_{2k}(M)\to\sk_{2k}(N)$ for each $d+1\leq k\leq n$.
\end{proposition}

\begin{proof} 
Suppose that $k=d+1$. Define $A$ by the composite 
\[A\colon\mathbb{Z}\langle f_{1},\ldots,f_{s}\rangle\overset{\cong}{\longrightarrow} 
     H_{2d+2}(M)\overset{\gamma^{\vee}}{\longrightarrow} H_{2d+2}(N)\overset{\cong}{\longrightarrow} 
     \mathbb{Z}\langle f'_{1},\ldots,f'_{s}\rangle\] 
where the left map is from~(\ref{lineariso}), $\gamma^{\vee}$ is the dual of $\gamma$, and the right map 
is the inverse of~(\ref{lineariso}) with respect to $N$. As $A$ is a composite of linear isomorphisms, it 
too is a linear isomorphism. Now consider the diagram
\[
\xymatrix{
\Z\langle{f_1,\ldots,f_s}\rangle\ar[r]^-{\cong}\ar[d]^-{A}	&H_{2d+2}(M)\ar[r]^-{g_M}\ar[d]^-{\gamma^{\vee}}	&\pi_{2d+1}(\sk_{2d}(M))\ar[d]^-{(\Phi_{2d})_*}\\
\Z\langle{f'_1,\ldots,f'_s}\rangle\ar[r]^-{\cong}			&H_{2d+2}(N)\ar[r]^-{g_N}								&\pi_{2d+1}(\sk_{2d}(N)).
}
\]  
The left square commutes by definition of $A$ and the right square commutes by~\eqref{eq_diag_MN}. The composites along the top and bottom rows are the 
definitions of $g$ and $g'$ respectively. By Proposition~\ref{Phi_2d}, $\Phi_{2d}$ is a homotopy 
equivalence. Thus the outer rectangle satisfies the hypotheses of 
Lemma~\ref{lemma_inductive argument}, implying that there is a homotopy equivalence 
\(\Phi_{2d+2}\colon\sk_{2d+2}(M)\to \sk_{2d+2}(N)\).

For $k>d+1$,  assume inductively that there is a homotopy equivalence 
\(\Phi_{2k}\colon\sk_{2k}(M)\to\sk_{2k}(N)\). 
Consider the diagram
\begin{equation}\label{Adiagram} 
\begin{split}
\xymatrix{
\Z\langle{f_1,\ldots,f_{s_k}}\rangle\ar[r]^-{\cong}\ar@{..>}[d]^-{A}	&H_{2k+2}(M)\ar[r]^-{g_M}  &\pi_{2k+1}(\sk_{2k}(M))\ar[d]^-{(\Phi_{2k})_*}\\
\Z\langle{f'_1,\ldots,f'_{s_k}}\rangle\ar[r]^-{\cong} 	&H_{2k+2}(N)\ar[r]^-{g_N} &\pi_{2k+1}(\sk_{2k}(N))
}\end{split}
\end{equation}  
where $A$ will be defined momentarily. Localize away from $\mathcal{P}$. Then 
the maps $g_M$ and $g_N$ are isomorphisms by Corollary~\ref{cor_g_M and g_BT}. The top 
and bottom rows are the definitions of $g$ and $g'$ as in \eqref{eq_map_g} and \eqref{eq_map_g_prime} respectively, so they are both isomorphisms. 
Define $A$ by the composite 
\[A\colon\mathbb{Z}\langle f_{1},\ldots, f_{s_k}\rangle\overset{g}{\longrightarrow} 
      \pi_{2k+1}(\sk_{2k}(M))\xrightarrow{(\Phi_{2k})_{\ast}}
      \pi_{2k+1}(\sk_{2k}(N))\xrightarrow{(g')^{-1}}\mathbb{Z}\langle f'_{1},\ldots,f'_{s_k}\rangle.\] 
Then $A$ is a linear isomorphism and it makes~(\ref{Adiagram}) commute. Thus~(\ref{Adiagram}) 
satisfies the hypotheses of the local version of Lemma~\ref{lemma_inductive argument} discussed 
in Remark~\ref{Aremark}, implying that there is a homotopy equivalence~\(\Phi_{2k+2}\colon\sk_{2k+2}(M)\to\sk_{2k+2}(N)\). 
This completes the induction.
\end{proof}

\newcommand{\etalchar}[1]{$^{#1}$}

\end{document}